\documentclass[12pt]{article}

\usepackage{amssymb}
\usepackage{amsthm}
\usepackage{amsmath}
\usepackage{verbatim}
\usepackage{mathrsfs}
\usepackage{enumerate}
\usepackage[usenames,dvipsnames] {color}

\usepackage{bm}

\let\cal=\mathcal      

\def\mcc{M\raise.5ex\hbox{c}C}
\def\mccarthy{M\raise.5ex\hbox{c}Carthy}

\def\ie{{\it i.e. }}

\def\h{{\cal H}}

\def\K{{\cal K}}

\def\m{Mult}

\def\l{\lambda}

\let\i=\infty

\def\la{\langle}
\def\ra{\rangle}
\def\={\ = \ }    
\def\ot{\otimes}

\def\A{{\cal A}}

\def\C{\mathbb C}

\def\dis{\displaystyle}

\def\be{\setcounter{equation}{\value{theorem}} \begin{equation}}
\def\ee{\end{equation} \addtocounter{theorem}{1}}
\def\beq{\begin{eqnarray*}}
\def\eeq{\end{eqnarray*}}

\def\bp{{\sc Proof: }}
\def\ep{{}{\hfill $\Box$} \vskip 5pt \par}

\def\bl{\begin{lemma}}
\def\el{\end{lemma}}
\def\bt{\begin{theorem}}
\def\et{\end{theorem}}
\def\bprop{\begin{prop}}
\def\eprop{\end{prop}}
\def\bd{\begin{definition}}
\def\ed{\end{definition}}
\def\br{\begin{remark}}
\def\er{\end{remark}}
\def\bexer{\begin{exercise}}
\def\eexer{\end{exercise}}
\def\bfig{\begin{figure}}
\def\efig{\end{figure}}

\numberwithin{equation}{section}

\title{Global holomorphic functions in several non-commuting variables II}
\author{Jim Agler
\thanks{Partially supported by National Science Foundation Grants
DMS 1361720 and DMS 1665260}
\\ U.C. San Diego\\ La Jolla, CA 92093
\and
John E. M\raise.5ex\hbox{c}Carthy
\thanks{Partially supported by National Science Foundation Grant  
DMS 1565243
}
\\ Washington University\\ St. Louis, MO 63130
}

\definecolor{dark_purple}{rgb}{0.4, 0.0, 0.4}

\def\bbm{\mathbb{M}}

\def\m{\mathcal{M}}
\def\n{\mathcal{N}}
\def\h{\mathcal{H}}

\def \b{\mathcal{B}}

\def\be{\begin{equation}}
\def\ee{\end{equation}}

\newcommand\de{\delta}
\renewcommand\bd{B_\delta}

\newcommand\md{\bbm^d}
\newcommand\hibd{H^\infty(B_\delta)}
\newcommand\hiobd{H_1^\infty(B_\delta)}
\def\pd{{\mathcal P}_d}
\newcommand\Vl{V_\lambda}
\newcommand\Il{I_\lambda}

\newcommand\lam{\lambda}

\def\s0{s_0}
\def\p0{p_0}

\DeclareMathOperator{\id}{id}
\DeclareMathOperator{\alg}{Alg}
\newcommand{\tensor}[2]{\text{ }{\begin{smallmatrix} #1 \\ \otimes\\ #2\end{smallmatrix}}\text{  }}
\newcommand{\triptensor}[3]{\text{ }{\begin{smallmatrix} #1 \\ \otimes\\ #2\\ \otimes\\ #3\end{smallmatrix}}\text{  }}

\renewcommand\O{\Omega}
\newcommand\mnd{{\bbm}_n^d}
\newcommand\mn{{\bbm}_n}
\newcommand\mm{{\mathbb M}}
\newcommand\cn{{{\mathbb C}^n}}
\renewcommand\L{{\mathcal L}}
\newcommand\holnc{{\rm Hol}^{\rm nc}_\h(\Omega)}

\begin{document}
\maketitle

%%%%%%%%%%%%%%%%%%%%%%%%%%%%%%%%

%%%%%%%%%%%%%%%%%%%%%%%%%%%%%%%%%%%%%%%%%%%%%%%%%%%%%%%%%%%%%%%%%%%%%%%%

\bibliographystyle{plain}
\allowdisplaybreaks

\theoremstyle{definition}
\newtheorem{defin}[equation]{Definition}
\newtheorem{lem}[equation]{Lemma}
\newtheorem{cor}[equation]{Corollary}
\newtheorem{prop}[equation]{Proposition}
\newtheorem{thm}[equation]{Theorem}
\newtheorem{theorem}[equation]{Theorem}
\newtheorem{claim}[equation]{Claim}
\newtheorem{ques}[equation]{Question}
\newtheorem{prob}[equation]{Problem}
\newtheorem{fact}[equation]{Fact}
\newtheorem{rem}[equation]{Remark}
\newtheorem{rems}[equation]{Remarks}
\newtheorem{notation}[equation]{Notation}
\newtheorem{exam}[equation]{Example}
\newtheorem{con}[equation]{Conjecture}

Abstract: We give a new proof that bounded non-commutative functions on polynomial polyhedra
can be represented by a realization formula, a generalization of the transfer function realization
formula for bounded analytic functions on the unit disk.

Math Subject Classification: 15A54

\section{Introduction}
\label{seca}

Let $\mn$ denote the $n$-by-$n$ matrices with complex entries, and let
$\md = \cup_{n=1}^\i \mnd$ be the set of all $d$-tuples of matrices of the same size.
An nc-function\footnote{nc is short for non-commutative} on a set $E \subseteq \md$ is a
function $\phi: E \to \mm^1$ that satisfies
\begin{enumerate}[(i)]
\item
$\phi$ is graded, which means that if $x \in E \cap \mnd$, then $\phi(x) \in \mn$.
\item
$\phi$ is intertwining preserving, which means if $x,y \in E$ and $S$ is a linear operator satisfying
$Sx = yS$, then $S \phi(x) = \phi(y) S$.
\end{enumerate}
The points $x$ and $y$ are $d$-tuples, so we write $x = (x^1, \dots, x^d)$ and
$y = (y^1, \dots, y^d)$. By $Sx = yS$ we mean that $S x^r = y^r S$ for each $1 \leq r \leq d$.
See \cite{kvv14} for a general reference to nc-functions.

The principal result of \cite{amfree} was a realization formula for nc-functions that are bounded on polynomial polyhedra;
the object of this note is to give a simpler proof of this formula, Theorem~\ref{thma1} below.

Let $\de$ be an $I$-by-$J$ matrix whose entries are non-commutative polynomials in $d$-variables. If $x \in \mnd$, then $\de(x)$ can be naturally thought of as an element of $\b(\C^J \otimes \C^n, \C^I \otimes \C^n)$, where $\b$ denotes the bounded linear operators, and all norms we use are operator norms on the appropriate spaces. We define
\be
\label{eqa1}
\bd \ := \ \{ x \in \md : \| \de(x) \| < 1 \} .
\ee
Any set of the form \eqref{eqa1} is called a \emph{polynomial polyhedron}.
Let $\hibd$ denote the nc-functions on $\bd$ that are bounded, and $\hiobd$ denote the closed unit ball,
those nc-functions that are bounded by $1$ for every $x \in \bd$.

\begin{defin}
A free realization for $\phi$ consists of an auxiliary Hilbert space $\m$ and  an isometry
\be
\label{eqa2}
\bordermatrix{~ &\C & \m \otimes \C^I \cr
\C&A&B\cr
\m \otimes \C^J & C&D
}
\ee
such that, for all $x \in \bd$, we have
\be
\label{eqa3}
\phi(x) \= \tensor{A}{1} + \tensor{B}{1} \tensor{1}{\de(x)} 
\left[ 1 - \tensor{D}{1} \tensor{1}{\de(x)} \right]^{-1} \tensor{C}{1} .
\ee
\end{defin}
The $1$'s need to be interpreted appropriately. If $x \in \mnd$, then
%$\phi(x)$ is an operator from $\tensor{\K}{\cn} $ to $\tensor{\L}{\cn}$, and 
 \eqref{eqa3}
means
\[
\phi(x) \= \tensor{A}{\id_\cn} + \tensor{B}{\id_\cn} \tensor{\id_\m}{\de(x)} 
\left[\triptensor{\id_\m}{\id_{\C^I}}{\id_\cn}  - \tensor{D}{\id_{\cn}} \tensor{\id_\m}{\de(x)} \right]^{-1} \tensor{C}{\id_\cn} .
\]
We adopt the convention of \cite{ptd13} and write tensors vertically to enhance legibility. The bottom-most entry
corresponds to the space on which $x$ originally acts; the top corresponds to the intrinsic part of the model on $\m$.

The following theorem was proved in \cite{amfree}; another proof appears in \cite{bmv16b}.
\bt
\label{thma1}
The function $\phi$ is in $\hiobd$ if and only if it has a free realization.
\et

It is a straightforward calculation that any function of the form \eqref{eqa3} is in $\hiobd$. We wish to prove the converse.
We shall use two other results, Theorems~\ref{thma2} and \ref{thma3} below.

 If $E \subset \md$, we let $E_n$ denote $E \cap \mnd$.
If $\K$ and $\L$ are Hilbert spaces, a $\b(\K,\L)$-valued nc function   on a set $E \subseteq \md$ is a function $\phi$ such that
\begin{enumerate}[(i)]
\item
$\phi$ is $\b(\K,\L)$ graded, which means if $x \in E_n$, then $\phi(x) \in \b(\K \otimes \cn, \L \otimes \cn)$.
\item
$\phi$ is intertwining preserving, which means if $x,y \in E$ and $S$ is a linear operator satisfying
$Sx = yS$, then 
\[
\tensor{\id_\L}{S} \phi(x) = \phi(y) \tensor{\id_\K}{S} 
. \]
\end{enumerate}

\begin{defin}
An nc-model for $\phi \in \hiobd$ consists of an auxiliary Hilbert space $\m$ and a $\b( \C, \m \otimes \C^J)$-valued
nc-function $u$ on $\bd$ such that, for all pairs $x,y \in \bd$ that are on the same level (\ie both in $\bd \cap \mnd$ for some $n$)
\be
\label{eqa4}
1 - \phi(y)^* \phi(x) 
\=
u(y)^* \left[ \tensor{1}{ 1 - \de(y)^* \de(x) } \right] u(x) .
\ee
\end{defin}
Again, the $1$'s have to be interpreted appropriately. If $x, y \in \bd \cap \mnd$, then 
\eqref{eqa4} means
\[
\id_{ \cn} - \phi(y)^* \phi(x) 
\=
u(y)^* \left[ \tensor{\id_\m}{ \id_{\C^J \ot \cn} - \de(y)^* \de(x) } \right] u(x) .
\]
%If $E \subseteq \md$, let $E^{[2]} $ denote the set of pairs $\{ (x,y) \}$ of elements of $E$ that are 
%on the same level (\ie both in the same $E_n$). Say that $\phi$ has an nc-model on $E$ if \eqref{eqa4} holds for all $(x,y) \in E^{[2]}$.
%
%We shall say $E \subseteq \md$ is an nc-set if it is closed with respect to direct sums, and that  it is closed with respect
%to unitary conjugation if 
%\[
%\forall_n \  \forall x \in E_n \ \forall u\ {\rm unitary\ in \ } \mn, \qquad u^* x u \in E .
%\]

\bt
\label{thma2} 
%\begin{enumerate}[(i)]
%\item
A graded function on $\bd$ has an nc-model if and only if it has a free realization.
%\item
%If $E \subseteq \bd$ is an nc-set that is closed with respect
%to unitary conjugation and $\psi$ is an nc-function on $E$ that has an nc-model on $E$, then 
%$\psi$ extends to a function in $\hiobd$ that has a free realization.
%\end{enumerate}
\et
Theorem~\ref{thma2} 
was proved in
\cite{amfree}, but a simpler proof is given by 
S. Balasubramanian in \cite{bal15}. Let us note for future reference that the functions $u$
in \eqref{eqa4} are locally bounded, and therefore holomorphic \cite[Thm. 4.6]{amfree}.
% (he does not state part (ii), but it follows from his method of proof).

The finite topology on $\md$ (also called the disjoint union topology)
is the topology in which a set  $\O$ is open if and only if 
 for every $n$,  
$\O_n$ is open in the Euclidean topology on $\mnd$. 
If $\h$ is a Hilbert space, and $\O$ is finitely open,
we shall let $\holnc$ denote the $\b(\C, \h)$ graded nc-functions on 
$\O$ that are holomorphic on each $\O_n$\footnote{
A function $u$ is holomorphic in this context if for each $n$, for each $x \in \Omega_n$,
for each $h \in \mnd$, the limit $\dis \lim_{t \to 0} \frac{1}{t} (u(x+th) - u(x) )$ exists.}
A sequence of functions $u^k$ on $\O$ is finitely locally uniformly bounded if for each point
$\l \in \O$, there is a finitely open neighborhood of $\l$ inside $\O$ on which the sequence is
uniformly bounded.

The following wandering Montel theorem is
proved in \cite{ammont}. If $u$ is in $\holnc$ and $V$ is a unitary operator on $\h$, define $V * u$ by
\[
\forall_n \quad (V*u) |_{\Omega_n} = \tensor{V}{\id_\cn} u |_{\Omega_n}.
\]
\bt
\label{thma3}
Let $\O$ be finitely open, $\h$ a Hilbert space, and $\{u^k\}$   a finitely locally uniformly bounded sequence in $\holnc$.
Then there exists a sequence $\{U^k\}$
of unitary operators on $\h$ such that $\{U^k * u^{k}\}$ has a  subsequence that converges
finitely locally uniformly to a function in ${\rm Hol}^{\rm nc}_\h(\bd)$.
\et

Let $\phi \in \hiobd$. We shall prove Theorem~\ref{thma1} in the following steps.
\begin{enumerate}[I]
\item
For every $z \in \bd$, show that $\phi(z)$ is in $ \alg(z)$, the unital algebra generated by the elements of $z$.
\item
Prove that for every finite set $F \subseteq \bd$, there is an nc-model for a function $\psi$ that agrees
with $\phi$ on $F$.
\item
Show that these nc-models have a cluster point that gives an nc-model for $\phi$.
\item
Use Theorem~\ref{thma2} to get a free realization for $\phi$.
\end{enumerate}

Remarks: 
\begin{enumerate}
\item
Step I is noted in \cite{amfree} as a corollary of Theorem~\ref{thma1}; proving it independently
allows us to streamline the proof of Theorem~\ref{thma1}.
\item
To prove Step II, we use one direction of \cite[Thm 1.3]{amfreePick15} that gives necessary and sufficient conditions to solve
a finite interpolation problem on $\bd$. The proof of necessity of this theorem used Theorem~\ref{thma1} above,
but for Step II we only need the sufficiency of the condition, and the proof of this in \cite{amfreePick15} did not use
Theorem~\ref{thma1}.
\item
All three known proofs of Theorem~\ref{thma1} start by proving a realization on finite sets, and then somehow taking a limit. In \cite{amfree}, this is done by considering partial nc-functions; in \cite{bmv16b}, it is done by using 
non-commutative kernels to get a compact set in which limit points must exist; in the current paper, we use the wandering Montel theorem.
\end{enumerate}

\section{Step I}
\label{secb}

Let $\{ e_j \}_{j=1}^n$ be the standard basis for $\C^n$. For $x$ in $\mn$ or $\mnd$, let $x^{(k)}$ denote the direct sum of $k$ copies of $x$. If $x \in \mnd$ and $s$ is invertible in $\mn$, then $s^{-1} x s$ denotes the $d$-tuple
$(s^{-1} x^1 s, \dots, s^{-1} x^d s)$.

\begin{lem}
\label{lemb1}
Let $z \in \mnd$, with $\|  z \| < 1$.
%Let $\A$ be the algebra generated by $z$, let $w \in \mn$, and 
Assume $w \notin \alg(z)$.
Then there is an invertible $s \in \mm_{n^2}$ such that $\| s^{-1} z^{(n)} s \| < 1$ and $\| s^{-1} w^{(n)} s \| > 1$.
\end{lem}
\bp
Let $\A = \alg(z)$.
Since $w \notin \A$, and $\A$ is finite dimensional and therefore closed,
the Hahn-Banach theorem says that there is a matrix $K \in \mn$ such that ${\rm tr}(aK) = 0 \  \forall\  a \in \A$ and
${\rm tr}(wK) \neq 0$.
Let $u \in \C^n \ot \C^n$ be the direct sum of the columns of $K$, and $v = e_1 \oplus e_2 \oplus \dots e_n$.
Then for any $b \in \mn$ we have
\[
{\rm tr}(bK) \= \la b^{(n)} u, v \ra .
\]
Let $ \A \ot \id$ denote $\{ a^{(n)} : a \in \A \}$.
We have 
 $\la a^{(n)} u, v \ra = 0\  \forall\  a \in \A$ and
$\la w^{(n)}  u, v \ra \neq 0$.

Let $\n = (\A \ot \id) u$. This is an $\A \ot \id $-invariant subspace, but it is not $w^{(n)}$ invariant
(since $v \perp \n$, but $v$ is not perpendicular to $w^{(n)} u$).
So decomposing $\C^n \ot \C^n$ as $\n \oplus \n^\perp$, every matrix in $\A \ot \id$ has $0$ in the $(2,1)$ entry,
and $w^{(n)}$ does not.

Let $s = \alpha I_{\n} + \beta I_{\n^\perp}$, with $\alpha > > \beta > 0$.
Then
\[
s^{-1} \begin{bmatrix}
A & B \\
C & D
\end{bmatrix}
s
\=
\begin{bmatrix}
A & \frac{\beta}{\alpha} B \\
\frac{\alpha}{\beta} C & D
\end{bmatrix} .
\]

If the ratio $\alpha/\beta$ is large enough, then for each of the $d$ matrices $z^r$, the corresponding
$s^{-1}( z^r \ot \id) s$ will have strict  contractions in the (1,1) and (2,2) slots, and each $(1,2)$ entry
will be small enough so that the whole thing is a contraction.

For $w$, however, as the $(2,1)$ entry is non-zero, the norm of $s^{-1} w^{(n)} s$ can be made arbitrarily large.
\ep

\begin{lem}
\label{lemb2}
Let $z \in \bd \cap \mnd$, and $w \in \mn$ not be in $\A := \alg(z)$.
Then there  is an invertible $s \in \mm_{n^2}$ such that
 $ s^{-1} z^{(n)} s \in \bd$ and 
$\| s^{-1} w^{(n)} s \| > 1$.
\end{lem}
\bp
As in the proof of Lemma~\ref{lemb1}, we can find an invariant subspace $\n$ for $\A \ot \id$ that is not $w$-invariant.
Decompose $\delta(z^{(n)})$ as a map from $(\n \otimes \C^J) \oplus (\n^\perp \otimes \C^J)$
into $(\n \otimes \C^I) \oplus (\n^\perp \otimes \C^I)$.
With $s$ as in Lemma~\ref{lemb1}, and $\alpha >> \beta >0$, and $P$ the projection from $\C^n \ot \C^n$ onto $\n$, we get 
\be
\label{eqb1}
\delta(s^{-1} z^{(n)} s)
\=
\begin{bmatrix}
\tensor{P}{\id}  \delta(z^{(n)}) \tensor{P}{\id} && \frac{\beta}{\alpha} \tensor{P}{\id} \delta(z^{(n)}) \tensor{P^\perp}{\id} \\ \\
0 & & \tensor{P^\perp}{\id}  \delta(z^{(n)}) \tensor{P^\perp}{\id} 
\end{bmatrix}.
\ee
The matrix is upper triangular because every entry of $\delta$ is a polynomial, and
$\n$ is $\A$-invariant.
For $\alpha/\beta$ large enough, every matrix of the form \eqref{eqb1} with $z \in \bd$ is a contraction, so
$ s^{-1} z^{(n)} s \in \bd$. 
But $s^{-1} w^{(n)} s$ will contain a non-zero entry multiplied by $\frac{\alpha}{\beta}$,
so we achieve the claim.
\ep
\begin{thm}
\label{thmb1}
If $\phi$ is in $\hibd$, then $\forall z \in \bd$, we have $\phi(z) \in {\rm Alg}(z)$.
\end{thm}
\bp
We can assume that $z \in \bd$ and that $\| \phi \| \leq 1$ on $\bd$.
Let $w= \phi(z)$. If $ w \notin {\rm Alg}(z)$, then
by Lemma \ref{lemb2}, there is an $s$ such that $s^{-1} z^{(n)} s \in \bd$
and
$\| \phi(s^{-1} z^{(n)} s )\| = \| s^{-1} w^{(n)} s \| > 1$, a contradiction.
\ep

Note that Theorem~\ref{thmb1} does not hold for all nc-functions. In
\cite{amif16} it is shown that there is a class of nc-functions, called fat functions, for which the
implicit function theorem holds, but Theorem~\ref{thmb1} fails.

\section{Step II}
\label{secc}

%The arguments here come from \cite{amfreePick15}, but are simplified by Theorem~\ref{thma3}.
%An nc-set is a subset of $\md$ that is closed with respect to direct sums.
%A graded function on an nc-set is intertwining preserving if and only it it preserves direct sums and similarities.

Let $F = \{ x_1, \dots, x_N \}$. Define $\l = x_1 \oplus \cdots \oplus x_N$, and define
$w = \phi(x_1) \oplus \cdots \oplus \phi(x_N)$.
As nc functions preserve direct sums (a consequence of being intertwining preserving)
we need to find a function $\psi$ in $\hiobd$ that has an nc-model, and satisfies $\psi(\l) = w$.

Let $\pd$ denote the nc polynomials in $d$ variables, and 
define
\[
\Il \= \{ q \in \pd\, : \, q(\l) = 0 \}.
\]
Let
\[
\Vl \= \{ x \in \md \, : \, q(x) = 0 \ {\rm whenever\ } q \in \Il \} .
\]
We need the following theorem from \cite{amfreePick15}:

\bt
\label{thmc1}
Let $\lam \in \bd \cap \mnd$ and $w \in \mn$.
There exists a function $\psi$ in the closed unit ball of $H^\infty(\bd)$ such that 
%$f(\l_i) = w_i$ for each $i$
$\psi(\lam) = w$
if 
(i) $w \in \alg(\l)$, so there exists $p \in \pd$ such that $p(\lam) = w$.

(ii) $\sup \{ \| p(x) \| : x \in \Vl \cap \bd \} \ \leq \ 1$.

Moreover, if the conditions are satisfied, $\psi$ can be chosen to have a free realization.
\et

Since $\phi(\l) = w$, by Theorem~\ref{thmb1}, there is a free polynomial $p$ 
such that $p(\l) = w$, so condition (i) is satisfied.
To see condition {\em (ii)}, note that for all $x \in \Vl \cap \bd$. we have $p(x) = \phi(x)$.
Indeed, by Theorem~\ref{thmb1}, there is a polynomial $q$ so that
$q(\l \oplus x) = \phi(\l \oplus x)$. Therefore $q(\l) = p(\l)$, so, since $x \in \Vl$,
we also have $q(x) = p(x)$, and hence $p(x) = \phi(x)$. 
But $\phi$ is in the unit ball of $\hiobd$,
so $\| \phi(x) \| \leq 1$ for every $x$ in $\bd$.

So we can apply Theorem~\ref{thmc1} to conclude that there is a function $\psi$ in $\hibd$ that has a free realization, and that
agrees with $\phi$ on the finite set $F$.

Remark: The converse of  Theorem~\ref{thmc1} is also true. Given Theorem~\ref{thmb1}, the converse
is almost immediate.

\section{Steps III and IV}
\label{secd}

Let $\Lambda = \{x_j \}_{j=1}^\infty$ be a countable dense set in $\bd$. For each $k$, let $F_k= \{ x_1, \dots, x_k\}$.
By Step II, there is a function $\psi^k \in \hiobd$ that has a free realization and agrees with $\phi$ on $F_k$.
By Theorem~\ref{thma2}, there exists a Hilbert space $\m^k$  and a
$\b(\C, \m^k \otimes \C^J)$ valued nc-function $u^k$ on $\bd$ so that, for all $n$, for all $x,y \in \bd \cap \mnd$, we have
\be
\label{eqd1}
1 - \psi^k(y)^* \psi^k(x) 
\=
u^k(y)^* \left[ \tensor{1}{ 1 - \de(y)^* \de(x) } \right] u^k(x) .
\ee
Embed each $\m^k$ in  a common Hilbert space $\h$.
Since the left-hand side of \eqref{eqd1} is bounded, it follows that $u^k$ are locally bounded, so we can apply Theorem~\ref{thma3} 
to find a sequence of unitaries $U^k$ so that, after passing to a subsequence,
$U^k * u^k$ converges to a function $v$ in $\holnc$.
We have therefore that 
\be
\label{eqd2}
1 - \phi(y)^* \phi(x) 
\=
v(y)^* \left[ \tensor{1}{ 1 - \de(y)^* \de(x) } \right] v(x) 
\ee
holds for all pairs $(x,y)$ that are both in $\Lambda \cap \mnd$ for any $n$, so by continuity,
we get that \eqref{eqd2} is an nc-model for $\phi$ on all $\bd$, completing Step III.

Finally, Step IV follows by applying Theorem~\ref{thma2}.

\section{Closing remarks}

One can modify the argument to get a realization formula for $\b(\K,\L)$-valued bounded nc-functions
on $\bd$, or to prove Leech theorems (also called Toeplitz-corona theorems---see \cite{Le14} and \cite{kr14}).
For finite-dimensional $\K$ and $\L$, this was done in \cite{amfree};
for infinite-dimensional $\K$ and $\L$ the formula was proved  in \cite{bmv16b}, using results from \cite{bmv16a}.

\bibliography{../../references}

\begin{thebibliography}{10}

\bibitem{ammont}
J.~Agler and J.E. M\raise.45ex\hbox{c}Carthy.
\newblock Wandering {M}ontel theorems for {H}ilbert space valued holomorphic
  functions.
\newblock To appear. arXiv:1706:05376.

\bibitem{amfree}
J.~Agler and J.E. M\raise.45ex\hbox{c}Carthy.
\newblock Global holomorphic functions in several non-commuting variables.
\newblock {\em Canad. J. Math.}, 67(2):241--285, 2015.

\bibitem{amfreePick15}
J.~Agler and J.E. M\raise.45ex\hbox{c}Carthy.
\newblock Pick interpolation for free holomorphic functions.
\newblock {\em Amer. J. Math.}, 137(6):1685--1701, 2015.

\bibitem{amif16}
Jim Agler and John~E. M\raise.45ex\hbox{c}Carthy.
\newblock The implicit function theorem and free algebraic sets.
\newblock {\em Trans. Amer. Math. Soc.}, 368(5):3157--3175, 2016.

\bibitem{bal15}
Sriram Balasubramanian.
\newblock Toeplitz corona and the {D}ouglas property for free functions.
\newblock {\em J. Math. Anal. Appl.}, 428(1):1--11, 2015.

\bibitem{bmv16b}
J.A. Ball, G.~Marx, and V.~Vinnikov.
\newblock Interpolation and transfer function realization for the
  non-commutative {Schur-Agler} class.
\newblock To appear. arXiv:1602.00762.

\bibitem{bmv16a}
Joseph~A. Ball, Gregory Marx, and Victor Vinnikov.
\newblock Noncommutative reproducing kernel {H}ilbert spaces.
\newblock {\em J. Funct. Anal.}, 271(7):1844--1920, 2016.

\bibitem{kr14}
M.~A. Kaashoek and J.~Rovnyak.
\newblock On the preceding paper by {R}. {B}. {L}eech.
\newblock {\em Integral Equations Operator Theory}, 78(1):75--77, 2014.

\bibitem{kvv14}
Dmitry~S. Kaliuzhnyi-Verbovetskyi and Victor Vinnikov.
\newblock {\em Foundations of free non-commutative function theory}.
\newblock AMS, Providence, 2014.

\bibitem{Le14}
Robert~B. Leech.
\newblock Factorization of analytic functions and operator inequalities.
\newblock {\em Integral Equations Operator Theory}, 78(1):71--73, 2014.

\bibitem{ptd13}
J.E. Pascoe and R.~Tully-Doyle.
\newblock Free {Pick} functions: representations, asymptotic behavior and
  matrix monotonicity in several noncommuting variables.
\newblock arXiv:1309.1791.

\end{thebibliography}
\end{document}